\documentclass[12pt,reqno]{article}
\usepackage{amsmath,amsthm,amsfonts,amssymb,color}
\usepackage{mathrsfs,mathbbol,graphicx}
\newcommand{\scr}[1]{\mathscr #1}
\definecolor{wco}{rgb}{0.5,0.2,0.3}

\numberwithin{equation}{section} \theoremstyle{remark}

\title{{\bf A Newton-Fixed Point Homotopy Algorithm For Nonlinear Complementarity Problems With Generalized Monotonicity} \footnote{Email: yuncholjong@yahoo.com
} }
\author{{\bf Yunchol Jong$^{a,*}$, Wonil Kim$^{a}$ }\\
\footnotesize{$^a$ Center of Natural Science, University of Science, Pyongyang, DPR Korea}\\}

\begin{document}
\maketitle
\def\A{\mathcal A}
\def\R{\mathbb R}     \def\Z{\mathbb Z}           \def\ff{\frac}        \def\ss{\sqrt}
\def\N{\mathbb N}
\def\dd{\delta}       \def\DD{\Delta}               \def\rr{\rho}
\def\vrr{\varrho}
\def\vpp{\varphi}
\def\ve{\varepsilon}
\def\<{\langle}       \def\>{\rangle}             \def\GG{\Gamma}       \def\gg{\gamma}
\def\ll{\lambda}      \def\LL{\Lambda}           \def\nn{\nabla}        \def\pp{\partial}    \def\tt{\triangle}
\def\d{\text{\rm{d}}} \def\bb{\beta}             \def\aa{\alpha}        \def\D{\scr D}
\def\E{\scr E}        \def\si{\sigma}            \def\ess{\text{\rm{ess}}}                   \def\Ker{\text{\rm{Ker}}}
\def\beg{\begin}      \def\beq{\begin{equation}} \def\enq{\end{equation}}                    \def\F{\scr F}
\def\beq{\begin{subequations}} \def\enq{\end{subequations}}
\def\Hess{\text{\rm{Hess}}}                       \def\e{\text{\rm{e}}}  \def\m{\text{\rm{min}}} \def\M{\text{\rm{max}}}
\def\B{\scr B}        \def\ua{\underline a}       \def\OO{\Omega}        \def\b{\mathbf b}
\def\oo{\omega}       \def\tt{\tilde}            \def\Ric{\text{\rm{Ric}}}\def\Re{\text{\rm{Re}}}
\def\cut{\text{\rm{cut}}} \def\P{\mathbb P} \def\ifn{I_n(f^{\bigotimes n})}
\def\fff{f(x_1)\dots f(x_n)} \def\ifm{I_m(g^{\bigotimes m})}
\def\pm{\pi_{\mu}}   \def\p{\mathbf{p}}   \def\ml{\mathbf{L}}
 \def\C{\scr C}      \def\aaa{\mathbf{r}}     \def\r{r}
\def\gap{\text{\rm{gap}}} \def\prr{\pi_{\mu,\varrho}}  \def\r{\mathbf r}
\def\Z{\mathbb Z} \def\Sect{{\rm Sect}}
\def\ii{{\rm i}_{\pp M}} \def\Ent{{\rm Ent}} \def\EE{\mathbb E}
\def\L{\mathbf L}
\def\le{\leqslant}
\def\ge{\geqslant}

\def\dist{\text{\rm dist}}
\def\div{\text{\rm div}}

\def\Cap{\text{\rm Cap}}
\def\Cov{\text{\rm Cov}}
\def\cut{\text{\rm cut}}
\def\Dom{\text{\rm Dom}}
\def\gap{\mathbf{gap}}
\def\supp{\text{\rm supp}}
\def\Var{\text{Var}}

\def\sect{\text{\rm sect}}\def\H{\mathbb H}
\def\wt{\widetilde}\def\wh{\widehat}
\def\tcr{\textcolor{red}}
\begin{abstract}
  In this paper has been considered probability-one global convergence of NFPH (Newton-Fixed Point Homotopy) algorithm for system of nonlinear equations and has been proposed a probability-one homotopy algorithm to solve a regularized smoothing equation for NCP with generalized monotonicity. Our results provide a theoretical basis to develop a new computational method for nonlinear equation systems and complementarity problems. Some preliminary numerical experiments shows that our NFPH method is useful and promissing for difficult nonlinear problems.
\end{abstract}


 \vskip 1.5cm

\section{Introduction}
  NCP (Nonlinear Complementarity Problem) is to find  $x\in R^n$ such that
  \begin{equation}\label{1}
     x_i\geq0,\quad f_i(x)\geq0,\quad x_if_i(x)=0,\quad i=1,...,n.
  \end{equation}
    Zhao and Li\cite{ZL} studied several properties of a homotopy solution path associated with nonlinear quasi-monotone complementarity problems. They established a sufficient condition to assure the existence and boundedness of this homotopy solution path. Their results provide a theoretical basis to develop a new computational method for quasi-monotone complementarity problems. In \cite{BW}, Billups and Watson considered probability-one global convergence of an interior FPH (Fixed Point Homotopy) algorithm for bounded MCP(Mixed Compementarity Problem). Their idea is to reformulate the MCP as a system of equations using FB(Fischer-Burmeister)- NCP function and then solve smooth approximations of this system with a homotopy method. Billups \cite{B} has considered probability-one global convergence of  FPH algorithm using smoothing of FB-NCP function for MCP satisfying  a coercivity or generalized monotonicity and strict feasibility. Watson\cite{W} has considered probability-one global convergence of  FPH algorithm  for monotone complementarity problem. His method involved reformulating the NCP as a system of smooth equations and applying a homotopy method to solve this system. In the context of Newton-based methods, such smooth reformulations of complementarity problems are inferior to nonsmooth reformulations due to slow local convergence for degenerate solutions. In contrast, nonsmooth reformulations allow much faster (superlinear or quadratic) convergence to degenerate solutions. In \cite{HY}, Hotta and Yoshise considered global convergence of a non-interior homotopy algorithm using CHKS(Chen-Harker-Kanzow-Smale)-smoothing function in case that map $f$ is $P_0$-mapping and has an interior feasible point, or $f$ is monotone and has an interior feasible point.
   In this paper, we propose NFPH (Newton-Fixed Point Homotopy) algorithm for nonlinear system and consider its probability-one global convergence, and extend it to solve a regularized CHKS-smoothing of nonsmooth reformulation for NCP with generalized monotonicity.


\section{Homotopy method for nonlinear system }
\subsection{Homotopy map and global convergence}
\indent
Our probability-one homotopy algorithm is based on the following Theorem.
\beg{thm}\label{thm1}\emph{(\cite{BW} and \cite{B})} Let $F:R^n\rightarrow R^n$ be a  $C^2$-function and suppose there exists  $C^2$-map $\rho:R^m\times[0,1)\times R^n\rightarrow R^n$   such that

  (i)  the $n\times (m+1+n)$  Jacobian matrix $D\rho(a,\lambda,x)$  has rank $n$  on the set
$$ \rho^{-1}(0)=\{(a,\lambda,x)\in R^m\times[0,1)\times R^n|\rho(a,\lambda,x)=0\}$$

 (ii) for any fixed $a\in R^m$  and $\lambda=0$ , the equation $\rho_a(\lambda,x)\equiv \rho(a,\lambda,x) =0 $ has a unique  solution $x^a\in R^n$ ,

 (iii) $\rho(a,1,x)=F(x)$   for any fixed $a\in R^m$ ,

 (iv) $\rho^{-1}(0)$  is bounded for any fixed $a\in R^m$ .

   Then for almost all $a\in R^m$  (in the sense of Lebesgue measure) there exists a zero curve $\gamma_a$  of $\rho_a$ , along which the Jacobian matrix $D\rho_a$  has rank $n$ , emanating from $(0,x^a)$ and reaching a zero   $\bar{x}$ of $F$  at $\lambda=1$ . Moreover, $\gamma_a$  does not intersect itself and is disjoint from any other zeros of $\rho_a$ .
    \end{thm}
    The expression "reaching a zero" means that there exists a sequence of points  $\{(\lambda_k,x^k)\}$ in  $\gamma_a$, accumulating at $(1,\bar{x})$.
  The popular homotopy used often in practice is the FPH   defined by
  $$ \rho(a,\lambda,x)=\lambda F(x)+(1-\lambda)(x-a).$$
 In this paper, we consider NFPH defined by
\begin{equation}\label{2}
 \rho(a,\lambda,x)=\lambda F(x)+(1-\lambda)G(x,a),
\end{equation}
where $G(x,a)=F(x)-F(a)+A(x-a)$  and  $A$ is a symmetric and positive definite matrix.
In what follows all of consideration will be made under the following assumption:

\textbf{Assumption 1.} $F$ is a  $C^2$-mapping and there exists a symmetric and positive definite matrix $A$ such that  $F'(x)+A$  is nonsingular for every $x\in R^n$ .
 \beg{rem}\label{rem1} Under the assumption 1, the map  $\rho(a,\lambda,x)$  defined by  $(\ref{2})$ satisfies conditions (i)$\sim$(iii) of Theorem \ref{thm1}.
 \end{rem}
The following Theorem provides a sufficient condition for condition (iv) of Theorem \ref{thm1} to be satisfied.
\beg{thm}\label{thm2}Suppose that $F:R^n\rightarrow R^n$  satisfies the assumption 1 and there exists  $\tilde{x}\in R^n$ and  $M>0$ such that
\begin{equation}\label{3}
   \|\tilde{x}\|_{A^\frac{1}{2}}< M,
\end{equation}
where $A$  is such as in the assumption 1 and  $\|x\|_{A^\frac{1}{2}}=\sqrt{x^TAx}$.
  If it holds
\begin{equation}\label{e4}
    (x-\tilde{x})^TF(x)\geq0
\end{equation}
 for every $x\in R^n$ such that $\|\tilde{x}-x\|_{A^{\frac{1}{2}}}\geq2M$,  then each zero curve $\gamma_a$  of  $\rho_a(\lambda,x)\equiv \rho(a,\lambda,x)$ defined by (2.1) is bounded for every  $a\in R^n$ such that
$$a+A^{-1}F(a)\in B=\{x\in R^n|\|x\|_{A^\frac{1}{2}}<M\}$$
and there exists a zero curve $\gamma_a$ , emanating from $(0,a)$  and reaching a zero  $\bar{x}$ of $F$  at $\lambda=1$  for almost all $a\in R^n$  with $a+A^{-1}F(a)\in B$. In particular, if $F'(\bar{x})$  is nonsingular, $\gamma_a$  has finite arc length.
\beg{proof}Let $\|x\|_{A^\frac{1}{2}}\geq3M$ .  Then we have
$$ \|x-a'\|_{A^\frac{1}{2}}\geq\|x\|_{A^\frac{1}{2}}-\|a'\|_{A^\frac{1}{2}}\geq3M-M=2M$$
for every $a'\in B$ , and since $\|\tilde{x}-a'\|_{A^\frac{1}{2}}\leq \|\tilde{x}\|_{A^\frac{1}{2}}+\|a'\|_{A^\frac{1}{2}}<2M$  by $(\ref{3})$ , we have
$$(x-\tilde{x})^TA(x-a')=(x-a'+a'-\tilde{x})^TA(x-a')=$$
$$=\|x-a'\|_{A^\frac{1}{2}}^2-(\tilde{x}-a')^TA(x-a')\geq\|x-a'\|_{A^\frac{1}{2}}^2-\|\tilde{x}-a'\|_{A^\frac{1}{2}}\|x-a'\|_{A^\frac{1}{2}}=$$
$$=\|x-a'\|_{A^\frac{1}{2}}(\|x-a'\|_{A^\frac{1}{2}}-\|\tilde{x}-a'\|_{A^\frac{1}{2}})>0$$
using generalized Schwartz inequality. Consequently, we have
\begin{equation}\label{5}
  (x-\tilde{x})^TA(x-a')>0
\end{equation}
for every $a'\in B$  and  $x\in R^n$ with $\|x\|_{A^\frac{1}{2}}\geq3M$ . Thus, if $a'=a+A^{-1}F(a)\in B$ , we have   $(x-\tilde{x})^T(\lambda F(x)+(1-\lambda)G(x,a))>0$ for $x\in R^n$  with $ \|x\|_{A^\frac{1}{2}}\geq3M$, i.e.
\begin{equation}\label{6}
   (x-\tilde{x})^T\rho_a(\lambda,x)>0
\end{equation}
for any $\lambda\in [0,1]$ by $(\ref{e4})$ and $(\ref{5})$ because $\|x-\tilde{x}\|_{A^\frac{1}{2}}\geq2M$ . The inequality $(\ref{6})$ implies that $\rho_a(\lambda,x)\neq0$ for  $\lambda\in [0,1]$ and $x\in \{x\in R^n|\|x\|_{A^\frac{1}{2}}\geq3M\}$ . Therefore, zero curve $\gamma_a$  of $\rho_a$  is contained in set $[0,1]\times\{x\in R^n|\|x\|_{A^\frac{1}{2}}<3M\}$  and is bounded, which proves the first proposition of Theorem together with Theorem \ref{thm1}. If  $F'(\bar{x})$ is nonsingular, , which implies that $\gamma_a$  has finite arc length.
\end{proof}
\end{thm}
\beg{defn} A function  $F:R^n\rightarrow R^n$ is said to be pseudo-monotone at $\tilde{x}$  if $(x-\tilde{x})^T F(\tilde{x})\geq0$  for every $x$  implies that $(x-\tilde{x})^TF(x)\geq0$  for every $x$ .
\end{defn}
 \beg{cor}\label{cor1}  Suppose that the assumption 1 is satisfied and solution set of $F(x)=0$  is bounded, i.e. $\|\tilde{x}\|_{A^\frac{1}{2}}<M$   for every solution $\tilde{x}$ of $F(x)=0$ .  If $F(x)$  is pseudo-monotone at each solution $\tilde{x}$ , then there exists a zero curve $\gamma_a$ of $\rho_a$ , emanating from $(0,a)$ and reaching a zero $\bar{x}$  of $\rho_a$  at $\lambda=1$  for almost every $a\in R^n$ with $a+A^{-1}F(a)\in B$ , where   $A$ is such as in the assumption 1.
  \beg{proof} For every solution $\tilde{x}$  of  $F(x)=0$ and for every $x\in R^n$  with $\|x-\tilde{x}\|_{A^\frac{1}{2}}\geq2M$ , we have $(x-\tilde{x})^TF(\tilde{x})=0$ . Hence, by the pseudo-monotonicity of $F(x)$ at $\tilde{x}$  , we have
  \begin{equation}\label{7}
   (x-\tilde{x})^TF(x)\geq0
  \end{equation}
for every $x\in R^n$  with $\|x-\tilde{x}\|_{A^\frac{1}{2}}\geq2M$ , which implies $(\ref{e4})$. Thus, Theorem \ref{thm2} completes the proof of the corollary.
\end{proof}
\end{cor}
\beg{rem}\label{rem2} Suppose that $\tilde{x}$  is a solution of $F(x)=0$. If $F$  is monotone in $\tilde{x}$ , $F$  is pseudo-monotone at $\tilde{x}$ . If $F$  is pseudo-monotone at $\tilde{x}$, $F$  satisfies $(\ref{e4})$.
\end{rem}
\subsection{Tracking the zero curve}
The zero curve can be tracked by the procedure considered in section 3.2 of [2]. As discussed in section 2.1, the zero curve can be parameterized by arc length : Let $(\lambda(s),x(s))$ be the point on $\gamma_a$  of arc length $s$ away from $(0,a)$. Tracking the zero curve involves generating a sequence of points  $\{y^k\}\in R^{n+1}$, with $y^0=(0,a)$, that lie approximately on the curve in order of increasing arc length. That is, $y^k=(\lambda(s_k),x(s_k))$ , where $\{s_k\}$  is some increasing sequence of arc lengths. The subroutine STEPNX from HOMPACK90 [6] is used to handle the curve tracking. At each iteration, STEPNX uses a predictor-corrector algorithm to generate the next point on the curve. The prediction phase requires  the corresponding unit tangent vector to the curve,$(y')^k=(\lambda'(s_k),x'(s_k))$  for each iterate  $y^k$. This is accomplished by finding an element $\eta$ of the null space of $\nabla\rho_a(y^k)$  and setting $(y')^k=\eta/\|\eta\|$ or $(y')^k=-\eta/\|\eta\|$, where the sign is chosen so that $(y')^k$  makes an acute angle with $(y')^{k-1}$ , for $k > 0$. On the first iterate, the sign is chosen so that the first component(corresponding to $\lambda$) of  $(y')^0$ is positive.
At each iteration after the first, STEPNX approximates the zero curve with a Hermite cubic polynomial $c^k(s)$ , which is constructed using the last two points $y^{k-1}$  and $y^k$ , along with the associated unit tangent vectors $(y')^{k-1}$  and $(y')^k $. A step of length $h$ along this cubic yields the predicted point $w^{k,0}=c^k(s+h)$ . The first iteration uses a linear predictor instead, which is constructed using the starting point $y^0$  and its associated unit tangent vector.

Once the predicted point is calculated, a normal flow corrector algorithm [6] is used to return to the zero curve. Starting with the initial point $w^{k,0}$ , the corrector iterates $w^{k,j},j=1,2,...$ , are calculated via the formula $w^{k,j+1}=w^{k,j}+z^{k,j},j=0,1,...$ , where the step $z^{k,j}$  is the unique minimum-norm .

The corrector algorithm terminates when one of the following conditions is satisfied: the normalized correction step $z^{k,j}/(1+\|w^{k,j}\|)$  is sufficiently small, some maximum number of iterations (6 in our experiments) is exceeded, or a rank-deficient Jacobian matrix is encountered in the Newton equation. In the first case, set   $y^{k+1}=w^{k,j}$, calculate an optimal step size $h$ for the next iteration, and proceed to the next prediction step. In the second case, discard the point and return to the prediction phase, using a smaller step size if possible; otherwise, terminate curve tracking with an error return. In the third case, terminate the curve tracking, since $rank\nabla \rho_a < n$ should theoretically not happen and indicates serious difficulty. The step size $h$ is also never reduced beyond relative machine precision.

Finally, to emphasize the robustness and effectiveness of the proposed homotopy method, consider the one-dimensional equation $f(x) = 0$, where
 $$f(x) = \arctan(100x)/\pi + \sin(5x/(x^2 + 0.2))/2 + 0.1x.$$
 solution to the Newton equation . This function has a unique root at $x = 0$. For algorithms that rely on descent of a merit function, this root is difficult to find because the global minimum of the merit function $\theta(x)=\frac{f(x)^2}{2}$  is in a very narrow valley (Fig.1). Nevertheless, the probability-one homotopy algorithm of [2] easily found the root, tracking the homotopy zero curve in 32 steps from a starting point $x^0=0.5$ .

Our Newton-fixed point homotopy method with $A=I$  in $G(x,a)$ of $(\ref{2})$  found the root in  2 iterations(6 steps) and 3 iterations(13 steps) from the same starting point with step size $h=0.63$ and $h=1$, respectively.  As a comparison, PATH version 4.0 [7] was used from the same starting point. After 449 iterations, PATH terminated at $x$ =0.24233, corresponding to a local minimum of $\theta$. This function $\theta(x)$, while artificial, is representative of merit functions encountered in applications such as protein folding, analog circuit simulation, and aircraft configuration design. This simple numerical experiment shows that our NFPH algorithm is a promising method.
\begin{figure}[h]
\begin{center}
  \includegraphics{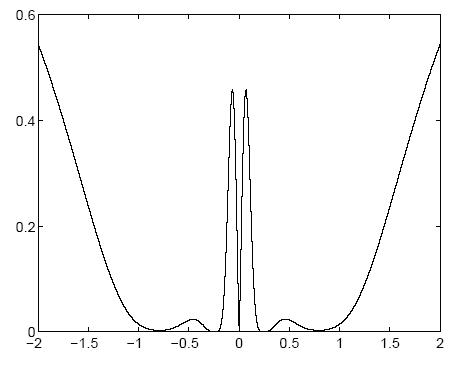}\\
  \caption{Graph of merit function $\theta(x) =\frac{f(x)^2}{2}$}\label{fig1}
\end{center}
\end{figure}

\section{Homotopy algorithm for NCP}
\beg{defn}
If a function $\varphi:R^2\rightarrow R$  is such that $\varphi(a,b)=0\Leftrightarrow a\geq0,b\geq0,ab=0$ , then   $\varphi$ is called NCP function.
\end{defn}
 The function $\varphi(a,b)=a+b-\sqrt{(a-b)^2}$ is a popular NCP function called min-function. Then   NCP (1.1) is equivalent to solve a system of nonlinear equations
$$\varphi\left(x_i,f_i(x)\right)=0,\quad i=1,...,n$$
and  the problem (1.1) is reduced to solve the following nonlinear system .
\begin{equation}\label{8}
    F(x,y)=\left(
             \begin{array}{c}
               f(x)-y \\
               \Phi(x,y) \\
             \end{array}
           \right),
\end{equation}
where $\Phi(x,y)=\left(\varphi(x_1,f_1(x)),...,\varphi(x_n,f_n(x))\right)^T$.  The min-function $\varphi(a,b)$  is nondifferentiable in case of  $a=b$  and this function is approximated by the following smooth function:
\begin{equation}\label{9}
    \varphi_\mu(a,b)=a+b-\sqrt{(a-b)^2+4\mu^2}
\end{equation}
 Let's approximate $(\ref{8})$ by the following smooth system using $(\ref{9})$:
 \begin{equation}\label{10}
    F^\mu(x,y)=\left(
             \begin{array}{c}
               f(x)-y+\mu x \\
               \Phi_\mu(x,y)+\mu y \\
             \end{array}
           \right),
 \end{equation}
where $\Phi_\mu(x,y)=\left(\varphi_\mu(x_1,f_1(x)),...,\varphi_\mu(x_n,f_n(x))\right)^T$ .  Then  $F^0(x,y)=F(x,y)$.
The smooth approximation often used for $(\ref{8})$ is
 $$
          \left(
             \begin{array}{c}
               f(x)-y \\
               \Phi_\mu(x,y) \\
             \end{array}
           \right).
 $$
 and $(\ref{10})$ is a regularized smoothing for $(\ref{8})$.
 \beg{lem}\label{lem1} For every $\mu>0,x$ and $y$, we have
 \begin{equation}\label{11}
   \|\Phi_\mu(x,y)-\Phi(x,y)\|\leq 2\mu\sqrt{n}.
 \end{equation}
 \beg{proof} When $\xi>0, \eta>0$, we have  $\sqrt{\xi+\eta}> \sqrt{\xi}$ and
$$
 ( \sqrt{\xi+\eta}- \sqrt{\xi})^2=\xi+\eta-2\sqrt{\xi+\eta}\sqrt{\xi}+\xi\leq 2\xi+\eta-2\sqrt{\xi}\sqrt{\xi}=\eta.
$$
Thus, letting $\xi_i=(x_i-y_i)^2,\quad \eta_i=(2\mu)^2,\quad i=1,...,n$ , it holds
$$\|\Phi_\mu(x,y)-\Phi(x,y)\|^2=\sum_{i=1}^n[\sqrt{(x_i-y_i)^2+4\mu^2}-\sqrt{(x_i-y_i)^2}\quad]^2\leq4n\mu^2 $$            which gives us the inequality $(\ref{11})$.

 \end{proof}
 \end{lem}
  Let's construct a homotopy to solve $(\ref{10})$ as follows.
\beg{equation}\label{12}
    \mu(\lambda)=\beta(1-\lambda),\quad \lambda\in [0,1],\quad \beta >0,
\end{equation}
\beg{equation}\label{13}
    \rho_a(\lambda,z)=\lambda F^{\mu(\lambda)}(z)+(1-\lambda)G^{\mu(\lambda)}(z,a),
\end{equation}
where $$z=\left(
            \begin{array}{c}
              x \\
              y \\
            \end{array}
          \right),\quad G^{\mu(\lambda)}(z,a)=F^{\mu(\lambda)}(z)-F^{\mu(\lambda)}(a)+A(z-a).$$

    Let's denote zero set of $\rho_a$  for fixed $a$  by
$$\rho_a^{-1}(0)=\{(\lambda,z)|\rho_a(\lambda,z)=0,\quad \lambda\in [0,1)\}$$
\beg{thm}\label{thm3} Suppose that the assumption 1 is satisfied for some matrix $A=cI$  and $F^{\mu(\lambda)}(z)$  defined by $(\ref{10})$, where $I$ is an identity matrix and $c$ is a positive constant.
  If $\|a\|_{A^\frac{1}{2}}< M,\quad \|\tilde{z}\|_{A^\frac{1}{2}}< M,\quad M\geq2\sqrt{cn}$  and
\begin{equation}\label{14}
 (z-\tilde{z})^TF(z)\geq0
\end{equation}
for every  $z\in R^{2n}$ with $\|z-\tilde{z}\|_{A^\frac{1}{2}}\geq2M$, then there exists a zero curve $\gamma_a$ of  $\rho_a$  defined by  $(\ref{13})$, emanating from $(0,a)$ and reaching a zero $\bar{z}$  of $F$  at $\lambda=1$ for almost all $a\in R^{2n}$  with $a+A^{-1}F^{\mu(\lambda)}(a)\in B=\{z\in R^{2n}|\|z\|_{A^\frac{1}{2}}< M\}.$
\beg {proof} First, let's prove that
\begin{equation}\label{15}
  (z-\tilde{z})^TF^{\mu(\lambda)}(z)\geq0
\end{equation}
for every $z\in R^{2n}$ with $\|z-\tilde{z}\|_{A^\frac{1}{2}}\geq2M$. Let $\mu$  denote $\mu(\lambda)$  for simplicity. We have
$$(z-\tilde{z})^TF^\mu(z)=(z-\tilde{z})^TF(z)+(z-\tilde{z})^T(F^\mu(z)-F(z))=$$
$$=(z-\tilde{z})^TF(z)+(x-\tilde{x})^T\mu x+(y-\tilde{y})^T(\Phi_\mu(z)-\Phi(z))+(y-\tilde{y})^T\mu y$$
$$=(z-\tilde{z})^TF(z)+\mu (z-\tilde{z})^Tz+(y-\tilde{y})^T(\Phi_\mu(z)-\Phi(z))$$
$$=(z-\tilde{z})^TF(z)+\mu\|z-\tilde{z}\|^2+\mu(z-\tilde{z})^T\tilde{z}+(y-\tilde{y})^T(\Phi_\mu(z)-\Phi(z))$$
$$\geq(z-\tilde{z})^TF(z)+\mu\|z-\tilde{z}\|^2-\mu\|z-\tilde{z}\|\|\tilde{z}\|-\|y-\tilde{y}\|\|\Phi_\mu(z)-\Phi(z)\|$$
$$\geq(z-\tilde{z})^TF(z)+\mu\|z-\tilde{z}\|^2-\mu\|z-\tilde{z}\|\|\tilde{z}\|-\|z-\tilde{z}\|\|\Phi_\mu(z)-\Phi(z)\|,$$
which gives us
\begin{equation}\label{16}
 \begin{split}
   (z-\tilde{z})^TF^\mu(z)& \geq(z-\tilde{z})^TF(z)+\mu\|z-\tilde{z}\|^2-\mu\|z-\tilde{z}\|\|\tilde{z}\|-2\mu\sqrt{n}\|z-\tilde{z}\| \\
     & =(z-\tilde{z})^TF(z)+\mu\|z-\tilde{z}\|(\|z-\tilde{z}\|-\|\tilde{z}\|-2\sqrt{n})
 \end{split}
\end{equation}
by $(\ref{11})$. It follows from $A=cI$  that $\|x\|_{A^\frac{1}{2}}=\sqrt{c}\|x\|$ , which implies $(\ref{15})$  by $(\ref{14})$  and $(\ref{16})$  because  $\|z-\tilde{z}\|\geq\frac{2M}{\sqrt{c}}\geq\|\tilde{z}\|+\frac{M}{\sqrt{c}}\geq\|\tilde{z}\|+2\sqrt{n}$. Therefore, a zero curve $\gamma_a$ of $\rho_a$  is bounded for each $a\in R^{2n}$  with $a+A^{-1}F^{\mu(\lambda)}(a)\in B=\{z\in R^{2n}|\|z\|_{A^\frac{1}{2}}< M\}$  by Theorem (2.2), and there exists a $\gamma_a$  of $\rho_a$ , emanating from $(0,a)$ and reaching a zero $\bar{z}$  of $F^{\mu(\lambda)}$  at $\lambda=1$ for almost every  $a\in R^{2n}$  with $a+A^{-1}F^{\mu(\lambda)}(a)\in B$.  Then, $\bar{z}$   is a zero of $F$  because  $F^{\mu(1)}(z)=F^0(z)=f(z)$  by $(\ref{10})$ and $(\ref{12})$.
\end{proof}
\end{thm}
\beg{lem}\label{lem2}\emph{([4])}. For every $\mu,\quad \varphi_\mu(a,b)=a+b-\sqrt{(a-b)^2+4\mu^2}=d$  if and only if $$\left(a-\frac{d}{2},b-\frac{d}{2}\right)\geq0,\quad \left(a-\frac{d}{2}\right)\left(b-\frac{d}{2}\right)=\mu^2.$$
\end{lem}
 Let $\lambda\in (0,1)$  and  $ A=\left(
                                 \begin{array}{cc}
                                   A_1 & 0 \\
                                   0 & A_2 \\
                                 \end{array}
                               \right) $.
 Then, it follows from $(\ref{13})$  that  for $a=\left(
                                                                           \begin{array}{c}
                                                                             a' \\
                                                                             a'' \\
                                                                           \end{array}
                                                                         \right),$

 $ \rho_a(\lambda,z)=0\quad\text{if and only if}$

  \begin{equation}\label{17}
    \Phi_\mu(z)=v,
  \end{equation}
  \begin{equation}\label{18.a}
    y=f(x)+r,\quad r=\mu x-(1-\lambda)[f(a')-a''+\mu a'-A_1(x-a')],
  \end{equation}
  \begin{equation}\label{18.b}
     v=-\mu y+(1-\lambda)[\Phi_\mu(a)+\mu a''-A_2(y-a'')]
  \end{equation}

 By Lemma \ref{lem2}, it follows from $(\ref{17})$ that $\rho_a(\lambda,z)=0$ if and only if
\begin{equation}\label{19}
\left(x_i-\frac{v_i}{2},y_i-\frac{v_i}{2}\right)\geq0,\quad \left(x_i-\frac{v_i}{2})(y_i-\frac{v_i}{2}\right)=\mu^2.
\end{equation}
\beg{lem}\label{lem3} Let $a=\left(
                            \begin{array}{c}
                              a' \\
                              a'' \\
                            \end{array}
                          \right)\geq\left(
                                       \begin{array}{c}
                                         \beta \\
                                         \beta \\
                                       \end{array}
                                     \right),$
where $\beta$  is the constant given in $(\ref{12})$. Then, we have $y>0$  for every $(\lambda,z)\in \rho^{-1}(0)$ with $z=(x,y)$.
\beg{proof} Let's prove that $\Phi_\mu(a)\geq0$ . If otherwise, there exists an $i$ such that $\varphi_\mu(a'_i,a''_i)=a'_i+a''_i-\sqrt{(a'_i-a''_i)^2+4\mu^2}<0$ .  Then, $a'_i+a''_i<\sqrt{(a'_i-a''_i)^2+4\mu^2}$  and it follows $a'_ia''_i<\mu^2.$ Thus, $\beta^2<\beta^2$  by $(\ref{12})$ and the condition of Lemma, which is a contradiction. Therefore, $\Phi_\mu(a)\geq0$ . Assume that there exists an index $i$  with $y_i\leq0$ . Then, it follows that $v_i>0$  by $(\ref{18.b})$ and $\Phi_\mu(a)\geq0$. But, by $(\ref{17})$,
$$v_i=\varphi_\mu(x_i,y_i)\leq\varphi(x_i,y_i)=2\min\{x_i,y_i\}\leq2y_i\leq0,$$
which is a contradiction. Thus, we have $y>0$ .
\end{proof}
\end{lem}
\beg{lem}\label{lem4} $a=\left(
                            \begin{array}{c}
                              a' \\
                              a'' \\
                            \end{array}
                          \right)\geq\left(
                                       \begin{array}{c}
                                         \beta \\
                                         \beta \\
                                       \end{array}
                                     \right),$
where $\beta$  is the constant given in $(\ref{12})$  and $(\bar{\lambda},\bar{z})$  be any limit of a smooth zero curve $\gamma_a\subset \rho^{-1}(0)$ . Let
$$P_1=\{i|\bar{x}_i=\infty\},\quad P_2=\{i|\bar{y}_i=\infty\},\text{and}\quad \{i|\bar{x}_i=-\infty\}.$$
  Then, $\bar{\lambda} <1$ implies that $i\ni P_2$   and  $f_i(\bar{x})=-\infty$ for every $i\in P_1$ and $P_2=N$. And $\bar{\lambda}=1$ implies that ( i ) $\bar{x}_i\leq0$  and $f_i(\bar{x})=\infty$  for every $i\in P_2$, (ii) $f_i(\bar{x})=\infty$  for every $i\in N$  and (iii) $\bar{y}_i=0$ and $f_i(\bar{x})\leq0$   for every $i\in P_1$ .
\beg{proof} 1) Consider the case of $\bar{\lambda}<1$. If $i\in P_1$, then $\bar{x}_i-\frac{\bar{v}_i}{2}=\infty$  where  $\bar{v}=\Phi_{\bar{\mu}}(\bar{z})$ and $\bar{\mu}=\mu(\bar{\lambda})$. Suppose the contrary. Then, there is a finite $C >0$ such that $\bar{x}_i-\frac{\bar{v}_i}{2}=C$ , and $\bar{v}_i=\infty$  because $\bar{x}_i=\frac{\bar{v}_i}{2}+C=\infty$ . Thus, $\bar{y}_i=-\infty$  by $(\ref{18.b})$. But, by Lemma 3.3, $\bar{y}_i\geq0$ which is a contradiction. Hence, by (3.13), $\bar{y}_i-\frac{\bar{v}_i}{2}=0$. Suppose that $i\in P_2$ . Then $\bar{y}_i=\infty$  and $\bar{v}_i=-\infty$  by $(\ref{18.b})$, which implies $\bar{y}_i-\frac{\bar{v}_i}{2}=\infty$, contradicting  $\bar{y}_i-\frac{\bar{v}_i}{2}=0$. Therefore, $i\ni P_2$ , i.e. $\bar{y}_i$  is a finite positive, from which it follows that $\bar{x}_i=\infty$  implies  $\bar{r}_i=\infty$, and $f_i(\bar{x})=-\infty$  by $(\ref{18.a})$. Now, let's prove $P_2=N$.  Let $i\in P_2$ .  Then, we have $\bar{v}_i=-\infty$  by $(\ref{18.b})$, and $\bar{y}_i-\frac{\bar{v}_i}{2}=\infty$, which implies $\bar{x}_i-\frac{\bar{v}_i}{2}=0$  and $\bar{x}_i=-\infty$  by $(\ref{19})$. Hence,   $i\in N$ and $P_2\subset N$ .  Let $i\in N$ . Then, $\bar{x}_i=-\infty$ and  $\bar{r}_i=-\infty$  by $(\ref{18.a})$. Because $\bar{x}_i\geq\frac{\bar{v}_i}{2}$  by $(\ref{19})$, we have $\bar{v}_i=-\infty$. Therefore, $\bar{y}_i=\infty$ , i.e. $i\in P_2$  by $(\ref{18.b})$, which implies  $N\subset P_2$. Thus, $P_2=N$.

2) Let's consider the case of $\bar{\lambda} =1$. In this case, by $(\ref{12})$ , we have $\bar{\mu} =0$, and by $(\ref{19})$
\begin{equation}\label{19'}
\left(\bar{x}_i-\frac{\bar{v}_i}{2},\bar{y}_i-\frac{\bar{v}_i}{2}\right)\geq0,\quad \left(\bar{x}_i-\frac{\bar{v}_i}{2})(\bar{y}_i-\frac{\bar{v}_i}{2}\right)=0.
\end{equation}
First, let's prove the proposition (i). If $i\in P_2$ , then  $\bar{v}_i\leq0$ and  $\bar{y}_i-\frac{\bar{v}_i}{2}=\infty$  by $(\ref{18.b})$, and
\begin{equation}\label{19''}
   \bar{x}_i-\frac{\bar{v}_i}{2}=0
\end{equation}
by $(\ref{19'})$. Consider the case of $\bar{v}_i=0$ . Then, $\bar{x}_i=0$  by $(\ref{19''})$, and $\bar{r}_i=0$  and  $\bar{y}_i=f_i(\bar{x})=\infty$  by $(\ref{18.a})$. Consider the case of $\bar{v}_i<0$ . Then, $\bar{x}_i<0$  by $(\ref{19''})$, and  $\bar{r}_i\leq0$ and $\infty=\bar{y}_i\leq f_i(\bar{x})$  by $(\ref{18.a})$. Second, let's prove the proposition (ii). If $i\in N$ , then $\bar{v}_i=-\infty$  by $(\ref{19'})$ and $\bar{y}_i=\infty$  by $(\ref{18.b})$, which implies that  $f_i(\bar{x})=\infty$  by the proposition (i). Finally, let's prove the proposition (iii). If $i\in P_1$ , we have  $ \bar{x}_i-\frac{\bar{v}_i}{2}=\infty$  by $(\ref{18.b})$ and
\begin{equation}\label{19'''}
 \bar{y}_i-\frac{\bar{v}_i}{2}=0
\end{equation}
 by $(\ref{19'})$. Suppose that  $i\in P_2$. Then, from the proposition (i), it follows that $\bar{x}_i\leq0$ , contradicting $i\in P_1$ . Hence, $i\ni P_2$  and $\bar{v}_i=0$  by $(\ref{18.b})$ because $\bar{\mu}=0$ and $\bar{\lambda}=1$.  Thus $\bar{y}_i=0$  by $(\ref{19'''})$, and we have $f_i(\bar{x})+\bar{r}_i =0$ and  $\bar{r}_i \geq0$  by $(\ref{18.a})$.  Consider the case of $\bar{r}_i =0$ . Then $\bar{y}_i =f_i(\bar{x})=0$.
 Consider the case of $\bar{r}_i >0$ . Then $f_i(\bar{x})+\bar{r}_i =0\geq f_i(\bar{x})$.
\end{proof}
\end{lem}
\beg{thm}\label{tm4} $a=\left(
                            \begin{array}{c}
                              a' \\
                              a'' \\
                            \end{array}
                          \right)\geq\left(
                                       \begin{array}{c}
                                         \beta \\
                                         \beta \\
                                       \end{array}
                                     \right),$
where $\beta$  is the constant given in $(\ref{12})$  and suppose that a map $f:R^n\rightarrow R^n$  satisfies   $f(a')>0$ and
\begin{equation}\label{20}
    (x-y)^T(f(x)-f(y))\geq0
\end{equation}
for some $\delta >0$  and for every  $x$ and $y$  with $\|x-y\|\geq\delta$ , and the assumption 1 is satisfied for  $F^{\mu(\lambda)}(z)$  defined by $(\ref{10})$.  Then, there is a smooth zero curve $\gamma_a\subset \rho^{-1}(0)$ , emanating from $(0,a)$  and reaching a NCP solution $\bar{z}=(\bar{x},\bar{y})$  with $\bar{y}=f(\bar{x})$  at $\lambda=1$ .
\beg{proof} If conditions of the Theorem are satisfied, conditions $(i)\sim(iii)$ of Theorem \ref{thm1} are satisfied by Remark \ref{rem1}. Then, there exists a smooth zero curve  $\gamma_a\subset \rho^{-1}(0)$ of $\rho_a(\lambda,z)$, along which the Jacobian matrix $D\rho_a(\lambda,z)$  has rank $n$ , emanating from  $(0,a)$ and not intersecting itself by Lemma 2.2 of [4].  Suppose that $\gamma_a$  is unbounded. Then there is an unbounded sequence $\{(\lambda_k,z^k)\}\subset \gamma_a$  such that $\lambda_k\rightarrow \tilde{\lambda}\in [0,1], z^k\rightarrow \tilde{z}$  and  $\|\tilde{z}\|=\infty$  as $k\rightarrow \infty$ , where $\tilde{z}=(\tilde{x},\tilde{y})$ . Let $\tilde{\lambda}<1$ . Then $\tilde{\lambda} =0$ is impossible because the zero curve  $\gamma_a$  can't return the starting point $(0,a)$  and the $\rho_a(0,z)$ has a unique simple zero $z=a$. Let $i\ni F=\{i|\tilde{z}_i \text{ is finite}\}$ . First, consider the case of $i\in P_1$ . Then, $f_i(\tilde{x})=-\infty$  by Lemma \ref{lem4}. Thus, we have
\begin{equation}\label{21}
(\tilde{x}_i-a'_i)(f_i(\tilde{x})-f_i(a'))=-\infty.
\end{equation}
Second, consider the case of $i\in N$ . Then, $\tilde{x}_i=-\infty$  and $\tilde{y}_i=\infty$  by Lemma \ref{lem4}, and it follows from $(\ref{18.a})$ that $f_i(\tilde{x})=\infty$ .  Therefore, we obtain the $(\ref{21})$  again. Third, consider the case when $\tilde{x}_i$  is finite. Then $\tilde{y}_i$   should be infinite because of $i\ni F$ , and $\tilde{y}_i=\infty$ . But, $P_2=N$  by Lemma \ref{lem4} and  $\tilde{x}_i=-\infty$, contradicting the finiteness of $\tilde{x}_i$ . Hence, we have $(\ref{21})$  for every $i\ni F$ .

 Suppose that $\|\tilde{x}-a'\|<\delta$  for every $\delta >0$. Then $\tilde{x}=a'$  because  $\delta$ is arbitrary. Hence, $\tilde{x}$  is bounded, which implies that $\tilde{y}$  should be unbounded. Since $P_2=N$ , if  $\tilde{y}_i=\infty$ for some $i$ , then  $\tilde{x}_i=-\infty$, contradicting  $\tilde{x}_i=a'_i$. Therefore, there is $\delta >0$ such that $\|\tilde{x}-a'\|\geq\delta$ . Then, by $(\ref{20})$, we have $(\tilde{x}-a')^T(f(\tilde{x})-f(a'))\geq0$ . Thus, taking account of $(\ref{21})$, there is $j\in F$  such that
 $$(\tilde{x}_i-a'_i)(f_i(\tilde{x})-f_i(a'))=\infty.$$
  The following two cases are possible.
  \begin{equation}\label{22}
    (i)\qquad \tilde{x}_j>a'_j,\quad f_j(\tilde{x})=\infty,
  \end{equation}
  \beg{equation}\label{23}
    (ii)\qquad \tilde{x}_j<a'_j,\quad f_j(\tilde{x})=-\infty.
  \end{equation}
   Consider the case (i). Because of $j\in F$, $\tilde{x}_j$  is bounded and so is $\tilde{r}_j$  by $(\ref{18.a})$. Then, we have $\tilde{y}_j=\infty$  by $(\ref{22})$ because $\tilde{y}_j=f_j(\tilde{x})+\tilde{r}_j$, contradicting $j\in F$.
   Consider the case (ii). Taking account of $(\ref{23})$, we have $\tilde{y}_j=-\infty$ similarly to the proof of the case (i), while $\tilde{y}_j\geq0$  by Lemma \ref{lem3}.  We have a contradiction again. Therefore, we have $\tilde{\lambda}=1$  and the equality $(\ref{21})$ for $i\in P_2\cup N$  by Lemma \ref{lem4}. For $i\in P_1$ , we have the equality $(\ref{21})$ again because $f(a')>0$ . Hence, it holds $(\ref{21})$ for every $i\ni F$ . Suppose that  $\|\tilde{x}-a'\|<\delta$  for every  $\delta>0$. Then $\tilde{x}=a'$ , which implies that $\tilde{y}$  is unbounded. Thus, $P_2\neq \varnothing $ and then $\tilde{x}_i\leq0$ for $i\in P_2$  by Lemma \ref{lem4}, contradicting $\tilde{x}_i=a'_i >0$. Therefore, there is $\delta >0$ such that $\|\tilde{x}-a'\|\geq\delta$ , and we have a contradiction again like the case of the $\tilde{\lambda}<1$ . Therefore, the zero curve $\gamma_a$  is bounded and it's accumulation point has form of $(1,\bar{z})$  by Lemma 2.3 of [4].  Then, $\bar{x}$  is a solution of $(\ref{1})$ and $\bar{y}=f(\bar{x})$ , where $\bar{z}=(\bar{x},\bar{y})$ . The proof is completed.

\end{proof}
\end{thm}
  The condition $(\ref{20})$ for map $f:R^n\rightarrow R^n$  is called generalized monotonicity.

\section{ Preliminary numerical experiments }
We can parameterize zero curve $\gamma_a$ of $\rho_a$  by its arclength $s$ , that is, there exist continuously differentiable functions $x(s)$ and $\lambda(s)$   such that
\beg{equation}\label{24}
    \varrho_a(x(s),\lambda(s))=0,
    x(0)=x^0,\lambda(0)=0.
\end{equation}

By differentiating the first equation of $(\ref{24})$, we obtain the following result: the homotopy path   is determined by the following initial value problem of the ordinary differential equations.
 \beg{equation}\label{26}
 \begin{array}{c}
   D\varrho_a(x(s),\lambda(s))\left(
                            \begin{array}{c}
                              \dot{x} \\
                              \dot{\lambda} \\
                            \end{array}
                          \right)=0,\\
   x(0)=x^0,\lambda(0)=0, \\
   \|\left(\dot{x}(s),\dot{\lambda}(s)\right)\|=1,\dot{\lambda}(0)>0,,
 \end{array}
\end{equation}
where $D\varrho_a(x(s),\lambda(s))=\left(\frac{\partial \rho_a}{\partial x},\frac{\partial \varrho_a}{\partial \lambda}\right).$                                                                                           Then the $x$ component of the solution $(x(s^*), \lambda(s^*))$ of $(\ref{26})$, which satisfies $\lambda(s^*) = 1$, is the solution which we need to find.
   We carried out numerical experiments for nonlinear equations to show the performance of our NFPH method by PC Pentium IV, 3.19GHz, 1.00GB of RAM.

   Our homotopy zero-curve trace was made by using ODE (ordinary differential equation) toolbox of MATLAB. In our experiment, we took $A=\alpha I, \alpha>0$ and made use of 'ode45' or 'ode15s', ODE solvers of MATLAB, to solve $(\ref{26})$.

Let $S_f$ be the upper bound used in the ODE part for arc length$s$ and
 $C_n$ be the number of intermediate checks $S_1, S_2, ...$ between 0 and $S_f$. Then interval $[0, S_f]$ is divided  in $(C_n +1)$ intervals of equal length. The zero-curve trace is finished as soon as one candidate solution has been found in any of the intervals defined by the checkpoints $[0, S_1, S_2, ... , S_f]$.
If path-following is successful, we obtain the more refined solution 'nsol' using 'fsolve', nonlinear equation solver of MATLAB, with initial guess 'hsol' obtained by homotopy method. In the following tables, $a$ is the starting point $x^0$, $sol$ is a solution of the given problem, $N_c$ is the number of checked intervals, and
$\quad fhom=F(hsol)/\left(1+\|hsol\|\right)$, $fnew=F(nsol)/\left(1+\|nsol\|\right)$.

Example 1.

$$F(x) = 2x - 4 + sin (2\pi x), x\in [-100, 100], a=0, sol=2.$$

Example 2.

$$F(x,q)=\left(
          \begin{array}{c}x^2 + q^2-1\\
           sin(x)-q)
           \end{array}
           \right),
    (x,q)\in [-100,100]\times[-100,100], a= [0,0]$$

Example 3.

$$F(x,q,z)= \left(
                            \begin{array}{c}
                              x + 0.5q + 0.3z-5 \\
                              0.6x + q + 0.1z-7 \\
                              0.2x + 0.4q + z-4
                            \end{array}
                          \right)$$
$$(x,q,z)\in [-100,100]\times[-100,100]\times[-100,100], a=[0,0,0]$$

Example4.
$$F(x) = arctan(100x)/\pi + sin(5x/(x2+0.2))/2 + 0.1x,
  x\in [-2,2], a=0.2, sol=0.$$
  The results of numerical experiments are shown in the following tables.
  $$ $$
   $\qquad Table 1(Example 1, S_f=2.5, C_n=70)$

\begin{small}
\begin{tabular}{|c|c|c|c|c|}
  \hline
  method & $N_c$ & hsol & nsol & time(s) \\
\hline
$\alpha=0.001$&	14&	2.0005&	2.00000& 1.0809\\
$\alpha=50$&	2&	1.9783&	2.00000&	0.3327\\
\hline
FPH & 20&	2.0008&	2.00000&	1.5547\\
\hline
NH & 14&	2.0000&	2.00000&	0.9153\\
\hline
\end{tabular}
\end{small}
$$ $$
$\qquad Table 2(Example 1, S_f=5, C_n=70)$

\begin{small}
\begin{tabular}{|c|c|c|c|c|}
  \hline
  method & $N_c$ & hsol & nsol & time(s) \\
\hline
$\alpha=0.001$&	7&	2.0005&	2.00000& 0.5911\\
$\alpha=50$&	1&	2.0032&	2.00000&	0.1862\\
\hline
FPH & 10&	2.0008&	2.00000&	0.7333\\
\hline
NH & 7&	2.0000&	2.00000&	0.5098\\
\hline
\end{tabular}
\end{small}
$$ $$
$\qquad Table 3(Example 2, S_f=20, C_n=70)$

\begin{small}
\begin{tabular}{|c|c|c|c|c|c|c|}
  \hline
   method & $N_c$ & hsol & nsol & fhom & fnew& time(s) \\
\hline
$\alpha=0.001$&	3
&	\begin{tabular}{c}
-7.4234e-001 \\
-6.7601e-001
 \end{tabular}
&\begin{tabular}{c}
	-7.3908e-001\\
-6.7361e-001
\end{tabular}
&\begin{tabular}{c}
	4.0222e-003\\
-2.7211e-006
\end{tabular}
&\begin{tabular}{c}
	3.5192e-011\\
3.6061e-012
\end{tabular}
& 0.5689\\
\hline
$\alpha=50$&1	
&\begin{tabular}{c}
-7.2803e-001\\
-7.3465e-001
\end{tabular}
&\begin{tabular}{c}
-7.3908e-001\\
-6.7361e-001
\end{tabular}
&\begin{tabular}{c}
3.4288e-003\\
3.4037e-003
\end{tabular}
&\begin{tabular}{c}
1.6008e-012\\
1.7203e-013
\end{tabular}
&1.4577\\
\hline
FPH&6	
&\begin{tabular}{c}
-7.3885e-001\\
-7.3465e-001
\end{tabular}
&\begin{tabular}{c}
-7.3908e-001\\
-6.7361e-001
\end{tabular}
&\begin{tabular}{c}
1.2271e-003\\
1.1222e-003
\end{tabular}
&\begin{tabular}{c}
2.3637e-012\\
5.1031e-013
\end{tabular}
&1.2431\\
\hline
NH &3	
&\begin{tabular}{c}
-7.3571e-001\\
-6.7112e-001
\end{tabular}
&\begin{tabular}{c}
-7.3908e-001\\
-6.7361e-001
\end{tabular}
&\begin{tabular}{c}
-4.1720e-003\\
-4.5865e-003
\end{tabular}
&\begin{tabular}{c}
4.1332e-011\\
4.2909e-012
\end{tabular}
&0.3311\\
\hline

\end{tabular}
\end{small}

$$ $$
$\qquad Table 4(Example 2, S_f=5, C_n=70)$

\begin{small}
\begin{tabular}{|c|c|c|c|c|c|c|}
  \hline
   method & $N_c$ & hsol & nsol & fhom & fnew& time(s) \\
\hline
$\alpha=0.001$&	10
&	\begin{tabular}{c}
-0.71376658\\
-0.65472611
 \end{tabular}
&\begin{tabular}{c}
-0.73908522\\
-0.67361213
\end{tabular}
&\begin{tabular}{c}
 -3.14294e -02\\
2.059234e-05
\end{tabular}
&\begin{tabular}{c}
1.3983767e- 07\\
1.4943515e-08
\end{tabular}
& 1.4112\\
\hline
$\alpha=50$&1	
&\begin{tabular}{c}
-0.72803942\\
-0.73465009
\end{tabular}
&\begin{tabular}{c}
-0.73908513\\
-0.67361202
\end{tabular}
&\begin{tabular}{c}
3.428826e-02\\
3.403779e-02
\end{tabular}
&\begin{tabular}{c}
1.6008305e-012\\
1.7202905e-013
\end{tabular}
&1.4576\\
\hline
FPH&21	
&\begin{tabular}{c}
-0.73915876\\
-0.67293244
\end{tabular}
&\begin{tabular}{c}
-0.739085321\\
-0.673612169
\end{tabular}
&\begin{tabular}{c}
-4.03203e-04\\
-3.67077e-04
\end{tabular}
&\begin{tabular}{c}
2.337165e-07\\
9.0814323e-010
\end{tabular}
&3.5987\\
\hline
NH &10	
&\begin{tabular}{c}
-0.7142857\\
-0.6550778
\end{tabular}
&\begin{tabular}{c}
-0.739085223\\
-0.673612122
\end{tabular}
&\begin{tabular}{c}
-3.080903e-02\\
-4.510372e-016
\end{tabular}
&\begin{tabular}{c}
1.288977e-07\\
1.378482e-08
\end{tabular}
&0.8984\\
\hline

\end{tabular}
\end{small}
$$ $$
$\qquad Table 5(Example 3, S_f=30, C_n=50)$

\begin{small}
\begin{tabular}{|c|c|c|c|c|c|c|}
  \hline
method & $N_c$ & hsol & nsol & fhom & fnew& time(s) \\
\hline
$\alpha=0.001$&	3
&	\begin{tabular}{c}
1.67897430 \\
 5.89119023\\
1.32550592
 \end{tabular}
&\begin{tabular}{c}
1.67155427 \\
 5.86510265\\
1.31964807
\end{tabular}
&\begin{tabular}{c}
3.0576e-003\\
4.2828e-003\\
2.4461e-003
\end{tabular}
&\begin{tabular}{c}
4.3102e-012\\
9.4505e-012\\
1.5067e-011
\end{tabular}
& 0.9735\\
\hline
$\alpha=50$&1	
&\begin{tabular}{c}
1.66456368\\
  5.92393350\\
1.31394188
\end{tabular}
&\begin{tabular}{c}
1.67155425\\
5.86510263\\
1.31964809
\end{tabular}
&\begin{tabular}{c}
2.8405e-003\\
7.4144e-003\\
2.2530e-003
\end{tabular}
&\begin{tabular}{c}
 -1.5842e-011\\
-3.3679e-011\\
-1.3984e-012
\end{tabular}
&0.9152\\
\hline
FPH&2	
&\begin{tabular}{c}
1.6731148\\
5.8574770\\
1.3209060
\end{tabular}
&\begin{tabular}{c}
1.6715542\\
5.8651026\\
1.3196482
\end{tabular}
&\begin{tabular}{c}
1.7504e-003\\
2.4505e-003\\
1.4003e-003
\end{tabular}
&\begin{tabular}{c}
-1.3416e-011\\
2.0760e-012\\
1.0126e-011
\end{tabular}
&0.6623\\
\hline
NH &3	
&\begin{tabular}{c}
1.6758000\\
5.8800000\\
1.3230000
\end{tabular}
&\begin{tabular}{c}
1.6715542\\
5.8651026\\
1.3196482
\end{tabular}
&\begin{tabular}{c}
1.7504e-003\\
2.4505e-003\\
1.4003e-003
\end{tabular}
&\begin{tabular}{c}
-1.3416e-011\\
2.0760e-012\\
1.0126e-011
\end{tabular}
&0.3325\\
\hline

\end{tabular}
\end{small}

$$ $$
$\qquad Table 6(Example 3, S_f=2, C_n=50)$

\begin{small}
\begin{tabular}{|c|c|c|c|c|c|c|}
\hline
method & $N_c$ & hsol & nsol & fhom & fnew& time(s) \\
\hline
$\alpha=0.001$&	33
&	\begin{tabular}{c}
1.50807356\\
5.29045112\\
1.19058642
 \end{tabular}
&\begin{tabular}{c}
1.67155426\\
5.86510262\\
1.31964809
\end{tabular}
&\begin{tabular}{c}
-0.07385090\\
  -0.10343821\\
  -0.05908049
\end{tabular}
&\begin{tabular}{c}
1.680643e-011\\
-1.5938702e-009\\
-6.454418e-010

\end{tabular}
& 6.7216\\
\hline
$\alpha=50$&1	
&\begin{tabular}{c}
1.69024906\\
5.79197983\\
1.33466574
\end{tabular}
&\begin{tabular}{c}
1.671554252\\
5.865102639\\
1.319648094
\end{tabular}
&\begin{tabular}{c}
-0.001861054\\
 -0.008413507\\
 -0.001461469
\end{tabular}
&\begin{tabular}{c}
6.8498442e-012\\
 2.679456e-012\\
5.0761474e-011
\end{tabular}
&0.7997\\
\hline
FPH&26	
&\begin{tabular}{c}
1.7520000\\
5.3961542\\
1.3848495
\end{tabular}
&\begin{tabular}{c}
1.67155424\\
5.86510263\\
1.31964809
\end{tabular}
&\begin{tabular}{c}
-0.01965900\\
  -0.06054967\\
  -0.01553925
\end{tabular}
&\begin{tabular}{c}
-4.948133e-010\\
 -4.022569e-010\\
 -3.449198e-010
\end{tabular}
&4.6629\\
\hline
NH &33	
&\begin{tabular}{c}
1.5047999\\
5.2800000\\
1.1880000
\end{tabular}
&\begin{tabular}{c}
1.67155425\\
5.86510264\\
1.31964809
\end{tabular}
&\begin{tabular}{c}
-0.075378066\\
 -0.105529292\\
 -0.060302452
\end{tabular}
&\begin{tabular}{c}
 6.0140810e-010\\
4.7703930e-010\\
6.6359281e-010
\end{tabular}
&1.6579\\
\hline

\end{tabular}
\end{small}
$$ $$
$\qquad Table 7(Example 4, S_f=5, C_n=70)$

\begin{small}
\begin{tabular}{|c|c|c|c|c|}
  \hline
  method & $N_c$ & hsol & nsol & time(s) \\
\hline
$\alpha=0.001$&	37&	-1.8652e-005&	3.1351e-011& 4.8317\\
$\alpha=1$&	18&	1.9012e-004&	7.6041e-020&	2.3264\\
$\alpha=75$&	1&	-2.3974e-004&	5.3332-022&	0.3118\\
\hline
FPH & 7&	-6.0612e-006  &	1.0732e-012&	0.8831\\
\hline
NH & 37&	5.0773e-005  &	2.9810e-022&	3.3849\\
\hline
\end{tabular}
\end{small}
$$ $$
As shown in Table 1 and Table 2, when $\alpha=50$, our $NFPH$ is much better than $FPH$ and $NH$ for Exampl 1.
For Example 2, our $NFPH$ with $\alpha=50$ is better than others in the accuracy of obtained solution (Table 3 and Table 4).
 Our $NFPH\quad(C_f=2,\quad \alpha=50)$ is much better than others in both of time and accuracy(Table 6) for Example 3. Especially, the $NFPH$ with $\alpha=75$ shows good performance for Example 4 (Table 7).
From the above numerical results, we can see that choice of proper $A$ can remarkably improve the performance of the $NFPH$.

\section{ Conclusion }

This paper describes a probability-one homotopy algorithm for solving nonlinear systems of equations and complementarity problems.  They are attractive because they are able to solve a qualitatively different class of problems than methods relying on merit functions. This claim is justified both theoretically and computationally. While the common homotopy used to solve nonlinear system $F(x)=0$  is FPH defined by $\rho(a,\lambda,x) = \lambda F(x)+(1-\lambda)(x-a)$, in this paper we considered probability-one global convergence of the algorithm based on NFPH defined by
$$\rho(a,\lambda,x) = \lambda F(x)+(1-\lambda)(F(x)-F(a)+A(x-a))$$
and extended the results to NCP with generalized monotonicity. The preliminary numerical experiments for some difficult nonlinear equations showed the robustness and fast convergence of the NFPH method. We expect the NFPH method would have better performance than NH(Newton Homotopy) or FPH method because NFPH combines the advantages of both NH and FPH.


 \beg{thebibliography}{99}

\bibitem{ZL} Y.-B. Zhao,  G.-N. Li,  \emph{Properties of a homotopy solution path for complementarity problems with quasi-monotone mappings, } Applied Mathematics and Computation 148, 93-104, 2004
\bibitem{BW} S.C. Billups, L.T. Watson. \emph{A probability-one homotopy algorithm for nonsmooth equations and mixed complementarity problems,} SIAM J. Optim., 12, 3, 606-626, 2002
\bibitem{B} S.C. Billups, \emph{A homotopy-based algorithm for mixed complementarity problems, } SIAM J. Optim., 12, 3, 583-605, 2002
\bibitem{W}  L.T. Watson, \emph{Theory of globally convergent probability-one homotopies for nonlinear programming,} SIAM J. Optim., 11, 3, 761-780, 2000
\bibitem{HY} K. Hotta, A. Yoshise, \emph{Global convergence of a class of noninterior point algorithms using Chen-Harker-Kanzow-Smale functions for nonlinear complementarity problems, } Math. Program., 86, 105-133, 1999
\bibitem{WSMMW} L. T. Watson, M. Sosonkina, R. C. Melville, A. P. Morgan, and H. F. Walker,  \emph{Algorithm
777: HOMPACK90: A suite of FORTRAN 90 codes for globally convergent homotopy algorithms,} ACM Trans. Math. Software, 23, 514-549, 1997
\bibitem{DF} S. P. Dirkse and M. C. Ferris, \emph{ The PATH solver: A non-monotone stabilization scheme for mixed complementarity problems,} Optim. Methods Softw., 5, 123-156, 1995

\end{thebibliography}

\end{document}